\documentclass[12pt]{amsart}

\usepackage{amsmath, amssymb, amsfonts, amsthm, amscd}
\usepackage{manfnt}
\usepackage{mathtools}
\usepackage{fullpage}
\usepackage[all]{xy}
   \SelectTips{cm}{10}
\usepackage{txfonts}
\newtheorem{thm}{Theorem}[section]

\newtheorem{prop}[thm]{Proposition}

\newtheorem{conj}[thm]{Conjecture}
\newtheorem{Scholium}[thm]{Scholium}

\newtheorem{prop-conj}[thm]{Proposition-Conjecture}

\theoremstyle{definition}
\newtheorem{defn}[thm]{Definition}

\theoremstyle{remark}
\newtheorem{rmk}[thm]{Remark}

\theoremstyle{remark}

\theoremstyle{remark}
\newtheorem{question}[thm]{Question}

\theoremstyle{remark}
\newtheorem{eg}[thm]{Example}

\newcommand{\Q}{\mathbb{Q}}
\newcommand{\Qb}{\overline{\mathbb{Q}}}
\newcommand{\Z}{\mathbb{Z}}
\newcommand{\CC}{\mathbb{C}}
\newcommand{\RR}{\mathbb{R}}
\newcommand{\bS}{\mathbb{S}}
\newcommand{\Ql}{\mathbb{Q}_{\ell}}
\newcommand{\Qlb}{\overline{\mathbb{Q}}_\ell}

\DeclareMathOperator{\Hom}{Hom}

\DeclareMathOperator{\Aut}{Aut}
\DeclareMathOperator{\Ind}{Ind}

\DeclareMathOperator{\Sym}{Sym}
\DeclareMathOperator{\Ad}{Ad}

\DeclareMathOperator{\Res}{Res}

\DeclareMathOperator{\Gr}{Gr}

\DeclareMathOperator{\pr}{pr}
\DeclareMathOperator{\Rep}{Rep}
\DeclareMathOperator{\Spec}{Spec}

\DeclareMathOperator{\Cent}{Cent}
\DeclareMathOperator{\Prim}{Prim}

\newcommand{\gal}[1]{\Gamma_{#1}} 
\newcommand{\Gal}{\mathrm{Gal}} 

\newcommand{\into}{\hookrightarrow}

\newcommand{\mc}{\mathcal}
\newcommand{\mf}{\mathfrak}
\newcommand{\mr}{\mathrm}

\newcommand{\tH}{\widetilde{H}}

\author{Stefan Patrikis}
\email{patrikis@math.mit.edu}
\address{Department of Mathematics\\MIT\\Cambridge, MA 02139}

\begin{document}
\title{Generalized Kuga-Satake theory and rigid local systems I: the middle convolution}
\thanks{I am very grateful to the organizers for the opportunity to speak at the conference. This work was conceived while unemployed in my parents' basement, but only brought to fruition with the support of NSF grant DMS-1303928}
\maketitle

\section{Background and motivation}
Let $X/\CC$ be a complex K3 surface. The classical construction of Kuga and Satake associates to $X$ a complex abelian variety $KS(X)/\CC$ along with a morphism
\[
H^2(X, \Q) \into H^1(KS(X), \Q) \otimes H^1(KS(X), \Q)
\]
of $\Q$-Hodge structures. Here is the Hodge-theoretic description of the construction, as rephrased by Deligne (\cite{deligne:weilK3}). For later purposes, we take $X$ to be projective, let $\eta$ be an ample line bundle on $X$, and replace the full cohomology group $H^2(X, \Q)$ by the $\eta$-primitive cohomology $\Prim^2_{\eta}(X, \Q)$ (the kernel of $\cup \eta$); for notational convenience we let $V_{\Q}$ be the weight-zero (Tate-twisted) $\Q$-Hodge structure 
\[
V_{\Q}= \Prim^2_{\eta}(X, \Q)(1).
\]
The Hodge bi-grading is then given by a representation of the Deligne torus $\mathbb{S}= \Res_{\CC/\RR} \mathbb{G}_m$ on $V_{\RR}$; this preserves the intersection pairing, and so we can set up the `lifting problem' of whether the dotted arrow in the following diagram can be filled in:
\begin{equation}\label{liftdiagram}
\xymatrix{
& \mr{GSpin}(V_{\RR}) \ar[d] \\
\bS \ar@{-->}[ur]^{\tilde{h}} \ar[r]^h & \mr{SO}(V_{\RR}).
}
\end{equation}
It is easily seen that this is possible, and that the ambiguity in such a lift is a choice of Tate twist. Now, associated to the quadratic space $V_{\Q}$ there is an (even--this is immaterial) Clifford algebra $C^+(V_{\Q})$, and there is a natural representation of $\mr{GSpin}(V_{\Q})$ on $C^+(V_{\Q})$. In these terms, the Kuga-Satake abelian variety is the complex abelian variety associated to the weight 1 polarizable $\Q$-Hodge structure whose underlying $\Q$-vector space is $C^+(V_{\Q})$, and whose Hodge bi-grading is given by the representation of $\bS$ induced by the lift (normalized to weight 1) $\tilde{h}$. 

It is a numerical miracle, crucially depending on the fact that $h^{2,0}(X)=1$, that the Hodge structure $(C^+(V_{\Q}), \tilde{h})$ has Hodge bi-grading of type $\{(1,0), (0, 1)\}$. A number of people have studied variants of the original construction (for example, \cite{morrison:KSabsurface}, \cite{galluzzi:KSabfourfold}, \cite{voisin:kuga-satake}, \cite{moonen:MTsurfaces}), always beginning with some Hodge structure resembling that of a K3 surface or abelian variety, and extracting some other Hodge structure (the `lift') that can by Riemann's theorem be shown to come from an abelian variety. 

But the diagram (\ref{liftdiagram}) above suggests a much broader question: abandon K3 surfaces and abelian varieties, and consider \textit{any} orthogonally polarized weight zero Hodge structure $V_{\Q}$ of motivic origin, for instance $\Prim^{2k}(X, \Q)(k)$ for \textit{any} smooth projective variety $X/\CC$. Again it is easy to find (unique up to Tate twist) a lift $\tilde{h}$, so we again have a (polarizable) $\Q$-Hodge structure $(C^+(V_{\Q}), \tilde{h})$. Is it motivic? Characterizing the essential image of the `Hodge realization' functor on the category $\mc{M}_{\CC}$ of pure homological motives over $\CC$ is as far as anyone knows a totally intractable problem; and indeed I know of no \textit{geometric} reason for believing $(C^+(V_{\Q}), \tilde{h})$ should come from algebraic geometry. If, however, the theory of motives is to reflect `l'identit\'{e} profonde entre la g\'{e}om\'{e}trie et l'arithm\'{e}tique,' we can hope to turn to arithmetic for guidance. 

Arithmetic in fact provides two compelling reasons for optimism, one Galois-theoretic and one automorphic. For a full discussion of these matters, see \cite{stp:variationsarxiv}; here I will review the Galois-theoretic aspect. Let $F$ be a number field, with an algebraic closure $\overline{F}$, and let $\Gamma_F= \Gal(\overline{F}/F)$ be its absolute Galois group. In parallel to the Hodge-realization   
\[
H^*_B \colon \mc{M}_{\CC} \to \Q\mr{HS}^{pol}
\]
to polarizable $\Q$-Hodge structures, there are $\ell$-adic realizations
\[
H^*_{\ell} \colon \mc{M}_F \to \Rep_{\Ql}(\Gamma_F)
\]
to the category of continuous representations of $\Gamma_F$ on finite-dimensional $\Ql$-vector spaces. In contrast to the Hodge-Betti realization, here there is a remarkable conjecture, due to Fontaine-Mazur, which in combination with the Tate conjecture\footnote{In the original paper \cite{fontaine-mazur}, the authors are careful to formulate a version of their conjecture that does not depend on the Tate conjecture.} would characterize the essential image of $H^*_{\ell}$.
\begin{defn}\label{FMgeom}
A Galois representation $\rho \colon \gal{F} \to \mr{GL}_N(\Qlb)$ is said to be \textit{geometric}\footnote{Some authors prefer the terminology `pseudo-geometric'} if for almost all finite place $v$ of $F$, the restriction of $\rho$ to a decomposition group $\gal{F_v}$ at $v$ is unramified; and for all $v \vert \ell$, $\rho|_{\gal{F_v}}$ is de Rham\footnote{The original definition by Fontaine-Mazur of `geometric' replaces `de Rham' by the \textit{a priori} stronger condition `potentially semi-stable'; but by the theorem of Berger (\cite{berger:dr=pst}) the two notions are equivalent.} in the sense of Fontaine (see \cite{brinon-conrad:cmi} for a textbook-style introduction to the theory of Fontaine).
\end{defn}
It follows from the base-change theorems for $\ell$-adic cohomology and Faltings' $p$-adic de Rham comparison isomorphism (\cite{faltings:crystalline}) that for any smooth projective $X/F$, $H^*_{\ell}(X)$ is geometric; in fact, the same holds for any $X/F$ separated of finite-type. Combining the conjectures of Tate and Fontaine-Mazur, we can predict that the essential image of $H^*_{\ell}$ on the category of pure motives over $F$ (now, to be precise, taken with $\Qb$-coefficients, and with a fixed embedding $\Qb \into \Qlb$) is those $\gal{F}$-representations that are semi-simple and geometric.

Now we can re-draw diagram (\ref{liftdiagram}) in its $\ell$-adic incarnation, letting $V_{\ell}$ be the $\ell$-adic realization of an orthogonally polarized motive over $F$, and replacing the representation $h$ of the Deligne torus $\bS$ with the (enormously richer) representation of $\gal{F}$ on $V_{\ell}$. Here, again, it is clarifying to strip away unnecessary features of the original context: what turns out to be important about the surjection $\mr{GSpin} \to \mr{SO}$ is that its kernel is equal to a \textit{central torus}. We are naturally led, then, to the following question.\footnote{This question first arose, with no hint of its connection to the Kuga-Satake construction, in the paper \cite{conrad:dualGW} of Brian Conrad.} Let $\tH \to H$ be any surjection of linear algebraic groups over $\Qlb$ whose kernel is a central torus in $\tH$, and suppose we are given a geometric Galois representation, or even one arising from the $\ell$-adic realization of a motive over $F$, $\gal{F} \xrightarrow{\rho} H(\Qlb)$; by an $H(\Qlb)$-valued geometric representation, we simply mean a homomorphism $\rho$ whose composition with some (equivalently, any) faithful finite-dimensional representation $H \into \mr{GL}_N$ is geometric in the sense of Definition \ref{FMgeom}. Does there exist a geometric lift $\tilde{\rho}$ filling in the following diagram?
\begin{equation}
\xymatrix{
& \tH(\Qlb) \ar[d] \\
\gal{F} \ar[r]^{\rho} \ar@{-->}^{\tilde{\rho}}[ur] & H(\Qlb).
}
\end{equation}
A beautiful theorem of Tate (\cite[Theorem 4]{serre:DSsurvey}), relying on the full force of global class field theory, asserts that the Galois cohomology group $H^2(\gal{F}, \Q/\Z)$ is zero. This implies (see \cite[Lemma 5.3]{conrad:dualGW}) that some, not necessarily geometric, lift $\tilde{\rho}$ exists; and also that when $\Qlb$ is replaced by $\CC$ (i.e., the Galois representations in question have finite image) it is possible to find geometric lifts (of course, here there is no `geometry,' since the only non-zero Hodge numbers will be $h^{0, 0}$). In general, there are some subtleties when $F$ admits real embeddings, but the answer to this geometric lifting question is essentially `yes': see \cite[Theorem 3.2.10]{stp:variationsarxiv} for the positive result over totally imaginary number fields (and \cite[Proposition 5.5]{stp:parities} for a complete description of the obstructions when $F$ has a real embedding; the details of this will not concern us). 

In combination with the Fontaine-Mazur-Tate conjecture, we are therefore led to the following conjecture. Let $\mc{M}_{F, E}$ denote the category of pure motives over $F$ with coefficients in some finite extension $E/\Q$; to have at our disposal an unconditionally Tannakian category of motives, we take this to mean Andr\'{e}'s category of motives for motivated cycles (see \cite{andre:motivated}), but one could alternatively assume the Standard Conjectures and take $\mc{M}_{F, E}$ to be the category of Grothendieck motives for any homological (or numerical) equivalence. A choice of fiber functor-- say, Betti cohomology, after embedding $F \into \CC$-- defines by Tannakian theory a motivic Galois group $\mc{G}_{F, E}$, a pro-reductive group over $E$, such that $\mc{M}_{F, E}$ is equivalent to the category of algebraic representations of $\mc{G}_{F, E}$.
\begin{conj}[See \S 4.3 of \cite{stp:variationsarxiv}]\label{KSconj}
Let $\tH \to H$ be a surjection of linear algebraic $E$-groups whose kernel is equal to a central torus in $\tH$, and let
\[
\rho \colon \mc{G}_{F, E} \to H
\]
be a motivic Galois representation. Then if either $F$ is totally imaginary, or the `Hodge numbers' of $\rho$ satisfy the (necessary) parity condition of \cite[Proposition 5.5]{stp:parities}, then there exists a finite extension $E'/E$ and a lifting of motivic Galois representations
\[
\xymatrix{
& \tH_{E'} \ar[d] \\
\mc{G}_{F, E'} \ar[r]_{\rho \otimes_E E'} \ar[ur]^{\tilde{\rho}} & H_{E'}.
}
\]
\end{conj} 
This conjecture should be viewed as a sharp arithmetic refinement (not enlarging $F$) of a vast generalization of the Kuga-Satake construction. It is of course also tempting to speculate that the analogous conjecture holds when we replace the number field $F$ by $\CC$; certainly the classical Kuga-Satake construction provides evidence for this. Along these lines, Serre has asked (\cite[8.3]{serre:motivicgalois}; this reference unfortunately has no explanation for this speculation) whether for any algebraically closed field $k$ ($\overline{F}$ or $\CC$, here) the derived group $\mc{G}_k^{der}$ is simply-connected.

The geometric aspect of Conjecture \ref{KSconj} (replacing $F$ by $\overline{F}$) already seems to be far out of reach, and the arithmetic refinement would seem to require combining a geometric result (over $\overline{F}$) with the Tate conjecture.\footnote{See \cite[Theorem 1.1.12]{stp:variationsarxiv} for an example of how this works in the classical Kuga-Satake setting, where the Tate conjecture for abelian varieties is available (Faltings' theorem).} The aim of the present paper is to provide examples of (a somewhat weakened version of) Conjecture \ref{KSconj}, examples in which the motives in question do not lie in the Tannakian sub-category $\mc{AV}_F$ of $\mc{M}_F$ generated by abelian varieties and Artin motives. If one wants to avoid difficult problems related to the Tate (or in our case `motivated-Tate') conjecture, it is convenient to phrase a weakened version of the problem as follows:
\begin{question}\label{question}
Let $\tH \to H$ be a surjection of linear algebraic $E$-groups with kernel equal to a central torus, as in Conjecture \ref{KSconj}. Given a motivic Galois representation $\rho \colon \mc{G}_{F, E} \to H$, does there exist
\begin{enumerate}
\item a finite extension $E'/E$;
\item a finite extension $F'/F$;
\item for each embedding $E' \into \Qlb$, inducing a place $\lambda' \vert \ell$ of $E'$, a geometric lifting
\[
\xymatrix{
& \tH(E'_{\lambda'}) \ar[d] \\
\gal{F} \ar[r]_{(\rho \otimes E')_{\lambda'}} \ar@{-->}^{\tilde{\rho}_{\lambda'}}[ur] & H(E'_{\lambda'}); 
}
\] 
\item and a faithful finite-dimensional representation $r \colon \tH_{E'} \into \mr{GL}_{N, E'}$ such that 
\[
r \circ \tilde{\rho}_{\lambda'}|_{\gal{F'}} \colon \gal{F'} \to \mr{GL}_N(E'_{\lambda'})
\]
is isomorphic to the $\lambda'$-adic realization of a motive, i.e. an object of $\mc{M}_{F', E'}$?
\end{enumerate}
Alternatively, letting $F= \CC$, we can ask the same question about the Hodge-Betti realizations, finding a rank $N$ motive whose $E'$-Hodge structure lifts the one given by the Betti realization of $\rho$. 
\end{question}
The truly tantalizing cases of this question are, in the spirit of the Vancouver conference, those arising from the `non-classical' context in which the motives (or Hodge structures) do not lie in $\mc{AV}_F$. We will content ourselves in this paper with producing many examples in the case (avoiding rationality questions) $\tH= \mr{GSpin}_5 \to H= \mr{SO}_5$, with $r= r_{spin}$ the (4-dimensional) spin representation;\footnote{A less interesting case also considered-- see Theorem \ref{so6}-- is $\mr{GSpin}_6 \to \mr{SO}_6$, with $r=r_{spin}$ the sum of the two half-spin representations.} these examples provide the first non-trivial examples of `generalized Kuga-Satake theory' beyond the Tannakian subcategory $\mc{AV}_F$ of $\mc{M}_F$. We will work with the $\ell$-adic version of Question \ref{question}, and therefore our results as written will only apply to (certain) motives over number fields; but it would also be possible to work purely Hodge-theoretically and deduce consequences for (certain) motives defined over $\CC$ but not $\Qb$. We postpone until \S \ref{rigid} a discussion of the precise setting of our motivic lifting theorem, and now conclude this introduction by giving one family of examples.

Consider the well-known family (with parameter $t$) of Calabi-Yau hypersurfaces
\[
\mathbf{X}_t= \{X_0^6+X_1^6+X_2^6+X_3^6+X_4^6+X_5^6= 6t X_0X_1X_2X_3X_4X_5\} \into \mathbb{P}^5.
\]
We denote by $\mathbf{X} \to \mathbb{P}^1 \setminus \{\mu_6, \infty \}$ this family of smooth projective sextic four-folds, whose $t$-fiber is $\mathbf{X}_t$. Each fiber $\mathbf{X}_t$ carries an action of the finite group scheme 
\[
\Gamma= \{(\zeta_0, \zeta_1, \ldots \zeta_5) \in \mu_6^6: \prod_i \zeta_i =1\} \subset \mu_6^6,
\]
with of course the diagonal subgroup $(\zeta, \zeta, \ldots, \zeta)$ acting trivially. For each $t \in (\mathbb{P}^1-\{\mu_6, \infty\})(\Qb)$, we can form the object of $\Gamma$-invariants in primitive cohomology
\[
M_t= \Prim^4(\mathbf{X}_t)^{\Gamma}(2) \in \mc{M}_{\Q(t)},
\]
where $\Q(t)$ is just the field of definition of the point $t$, and where we have Tate-twisted to weight zero. For any embedding $\Q(t) \into \CC$, the associated Hodge numbers are 
\[
h^{2,-2}= h^{1, -1}= h^{0,0}=h^{-1, 1}= h^{-2,2}=1, 
\]
and from this it is easy to see that $M_t$ does not lie in $\mc{AV}_{\Q(t)}$, nor after base-change does $M_t |_{\mc{G}_{\CC}}$ lie in $\mc{AV}_{\CC}$. $M_t$ gives rise to an orthogonal motivic Galois representation $\rho_t \colon \mc{G}_{\Q(t)} \to O(M_t)$, where the orthogonal pairing is given by Poincar\'{e} duality. There is an at most degree $2$ extension of $\Q(t)$ such that for all $\ell$, the restriction of the $\ell$-adic realization $\rho_{t, \ell}$ factors through the special orthogonal group. Here, then, is a sample result:
\begin{thm}\label{dwork}
For all $t \in \mathbb{P}^1(\Qb)-\{\mu_6, \infty\}$, there exists a number field $F$, in fact a quadratic extension (depending on $t$) of $\Q(\zeta_{24})$ such that:
\begin{enumerate}
\item for all places $\lambda$ of $\Q(\zeta_{24})$, there is a lift $\tilde{\rho}_{t, \lambda}$ of $\rho_{t, \lambda}$:
\[
\xymatrix{
& \mr{GSpin}(M_{t, \lambda}) \ar[d] \\
\gal{F(t)} \ar[ur]^{\tilde{\rho}_{t, \lambda}} \ar[r]_{\rho_{t, \lambda}} & \mr{SO}(M_{t, \lambda});
}
\]
\item and an object $\widetilde{M}_t$ of $\mc{M}_{F(t), \Q(\zeta_{24})}$ such that the $\lambda$-adic realization $\widetilde{M}_{t, \lambda}$ is isomorphic (as $\gal{F(t)}$-representation) to $r_{spin} \circ \tilde{\rho}_{t, \lambda}$.
\end{enumerate}
\end{thm}
Informally, the motives $M_t$ each admits a generalized Kuga-Satake lift. Note that the field $F$ has bounded degree, independent of $t$ and $\lambda$. This is some consolation for its unfortunate occurence in the statement of the theorem; sometimes with more work sharper results are possible.
\section{Rigid local systems and Katz's middle convolution algorithm}\label{rigid}
What is special about each of the motives $M_t$ arising from the Dwork family of hypersurfaces, substituting for the fact that we no longer have at our disposal the Hodge theory of abelian varieties? Roughly speaking, it is that the $\ell$-adic realization $M_{t, \ell}$ is isomorphic to a `fiber' of a rigid local system over suitably punctured $\mathbb{P}^1$. We will now take some time to describe Katz's beautiful theory of such rigid local systems. Let $k$ be an algebraically closed field, and let $\ell$ be a prime not equal to the characteristic of $k$. Let $S \subset \mathbb{P}^1$ be a finite set of closed points, and let $j \colon U= \mathbb{P}^1-S \into \mathbb{P}^1$ be the open complement. For $s \in S$, we let $I_s$ denote the inertia group at $s$. To put our results in the proper context, we will introduce not just Katz's notions of rigidity for local systems, but also their (obvious) extensions to $H$-local systems for any linear algebraic group $H$ over $\Qlb$ (see too \cite{yun:iccmsurvey} for these definitions). Recall that an $H$-local system on $U$ is simply a continuous homomorphism $\rho \colon \pi_1(U, x) \to H(\Qlb)$, where $x$ is a fixed geometric point of $U$. Two $H$-local systems are isomorphic if they are conjugate by an element of $H(\Qlb)$. There are at least two relevant notions of rigidity for such an $H$-local system:
\begin{defn}
An $H$-local system $\rho$ on $U$ is $H$-physically rigid if for all $H$-local systems $\rho'$ on $U$ such that $\rho|_{I_s} \cong \rho'|_{I_s}$ for all $s \in S$, we in fact have $\rho \cong \rho'$.\footnote{We paraphrase this condition as saying that `if $\rho$ and $\rho'$ are everywhere-locally conjugate, then they are globally conjugate.'} We say that $\rho$ is $H$-cohomologically rigid if
\[
H^1(\mathbb{P}^1, j_* \Ad(\rho))=0,
\]
where of course $\Ad \colon H \to \mr{GL}(\mr{Lie}(H))$ denotes the adjoint representation.
\end{defn}
If $\rho$ is tamely ramified at all $s \in S$-- this is the only case that will concern us-- then the Euler-Poincar\'{e} formula simply says
\[
\chi(\mathbb{P}^1, j_* \Ad(\rho))= \dim H \cdot \chi(U) + \sum_{s \in S} \dim \Ad(\rho)^{I_s}.
\]
If $H$ is moreover reductive, then $\chi(\mathbb{P}^1, j_* \Ad(\rho))= 2\dim \Ad(\rho)^{\pi_1(U, x)} - \dim H^1(\mathbb{P}^1, j_* \Ad(\rho))$, so $\rho$ is cohomologically rigid if and only if
\begin{equation}\label{cohrig}
\dim H \cdot (2-|S|)+ \sum_{s \in S} \dim \Ad(\rho)^{I(s)}= 2 \dim \Ad(\rho)^{\pi_1(U, x)}.
\end{equation}
Note that the analogue of irreducibility (when $H= \mr{GL}_N$) is the condition $\Cent_H(\rho)= Z(H)$, in which case the right-hand side of this formula is just $2 \dim Z(H)$.

For $H= \mr{GL}_N$, Katz has shown (\cite[Theorem 5.0.2]{katz:rls}) that if an irreducible lisse sheaf $\mc{F}$ is cohomologically rigid, then it is physically rigid; and conversely that if $k \subset \CC$, and if the analytic sheaf $\mc{F}^{an}$ associated to $\mc{F}$ is physically rigid on $U^{an}$, then $\mc{F}$ is cohomologically rigid (note that the `passage to the analytic' functor is fully faithful but not essentially surjective-- see \cite[Proposition 5.9.2]{katz:rls} for a discussion-- so that physical rigidity of $\mc{F}^{an}$ does not follow immediately from physical rigidity of $\mc{F}$; for the proof of this converse, see \cite[Theorem 1.1.2]{katz:rls}). 
Note that the implication `cohomologically rigid implies physically rigid' does not continue to hold for general $H \neq \mr{GL}_N$; this has the potential to complicate significantly the kinds of arguments given in this paper, although here (essentially because of low-dimensional coincidences) it does not intervene.

Underlying the approach of this paper, and our hope to exploit rigid local systems for purposes of the Kuga-Satake problem, is a (variant of a) conjecture of Simpson (see \cite[p.9]{simpson:higgs}), which we will state as a guiding philosophy rather than as a precise conjecture: 
\begin{Scholium}\label{rigimpliesmot}
Let $H$ be a reductive group, and let $\rho$ be an $H$-rigid local system on $U$ with quasi-unipotent local monodromies. Then $\rho$ is `motivic,' i.e. arises as a direct factor of the monodromy representation on the cohomology of a family of varieties $Y \to U$. 
\end{Scholium}
Note that we have not specified which kind of rigidity. For the purposes of heuristic argument, let us blur the two. Extrapolating from this principle, if we have a (known to be motivic) Galois representation or Hodge structure $\rho_t$ isomorphic to the $t$-fiber of some $H$-rigid local system $\rho$ on $U$, and we have a surjection $\tH \to H$, then to find a Kuga-Satake lift we should lift $\rho$ to $\tH$ (this is obviously possible over algebraically closed $k$; in the Galois case, when we work over a number field, we achieve this after a finite, but \textit{a priori} hard to specify, base-change). The lift $\tilde{\rho}$ will be $\tH$-rigid (compare equation \ref{cohrig} for $H$ and $\tH$), so it ought to be motivic. The (motivic) $t$-fiber $\tilde{\rho}_t$ should then provide the Kuga-Satake lift of $\rho_t$.

Of course, nothing so general as Scholium \ref{rigimpliesmot} is proven; but, remarkably, Katz has proven it, for irreducible $\rho$, in the case $H= \mr{GL}_N$! Bogner and Reiter (\cite{bogner-reiter:rigidsp4}) have combined Katz's work with a clever trick to establish a similar result for $H= \mr{Sp}_4$ as well. These two results are essentially all we will require. For the rest of this section, we proceed more formally and explain Katz's work as needed for our purposes.

As before, let $k$ be an algebraically closed field, and let $\ell$ be a prime invertible in $k$. We say that a lisse $\Qlb$-sheaf $\mc{F}$ on $U/k$ has index of quasi-unipotence dividing $N$ if the local monodromies, i.e. for each $s \in S$ the action of a choice of topological generator of $I_s \cong \hat{\Z}$, are quasi-unipotent with eigenvalues contained in the $N^{th}$ roots of unity. We restrict to $N$ invertible in $k$. Rigid local systems over $k$ will admit arithmetic descents, and for this purpose we consider the ring $R_{N, \ell}= \Z[\zeta_N, \frac{1}{N\ell}]$, with a \textit{fixed} choice of embedding $R_{N, \ell} \into k$. Now, we may always apply an automorphism of $\mathbb{P}^1$ to assume that $\infty \in S$, and we choose an ordering of the remaining points of $S$, writing $U= \mathbb{A}^1-\{\alpha_1, \ldots, \alpha_n\}$.  It turns out that a rigid local system on $U$ extends to a `universal' local system in which the points $\{\alpha_i\}$ are allowed to move in the plane, so we introduce the arithmetic configuration space
\[
S_{N, n, \ell}= R_{N, \ell}[T_1, \ldots, T_n][\frac{1}{\prod_{i \neq j} (T_i - T_j)}].
\]
The `universal' version of a given rigid local system will then live on the relative affine line with $n$ sections deleted, 
\[
\mathbb{A}^1_{S_{N, n, \ell}}-\{T_1, \ldots, T_n\}= \Spec R_{N, n, \ell}[T_1, \ldots, T_n, X]\left[ \frac{1}{\prod_{i \neq j}(T_i-T_j)} \cdot \frac{1}{\prod_{i} (X-T_i)} \right].
\]
Specialization to a local system on $U/k$ is then achieved by extending the fixed $R_{N, \ell} \into k$ to 
\begin{align*}
\phi \colon S_{N, n, \ell} &\to k \\
T_i &\mapsto \alpha_i.
\end{align*}
Now we can state (a slightly less precise version of) Katz's theorem on the `motivic' description of rigid local systems:
\begin{thm}[Theorem 8.4.1 of \cite{katz:rls}]\label{katzmotivic}
Let $\mc{F}_k$ be a lisse $\Qlb$-sheaf on $U/k$ that is 
\begin{itemize}
\item tamely ramified;
\item irreducible; 
\item cohomologically rigid; 
\item quasi-unipotent of index dividing $N$.
\end{itemize}
Fix a faithful character $\chi \colon \mu_N(R_{N, \ell}) \into \Qlb^\times$ (equivalently, an embedding $R_{N, \ell} \into \Qlb$). Then there exists a smooth affine family of $R_{N, \ell}$-schemes
\[
\pi \colon \mr{Hyp} \to \mathbb{A}^1_{S_{N, n, \ell}}-\{T_1, \ldots, T_n\}
\]
such that
\begin{itemize}
\item the fibers of $\pi$ are geometrically connected of some dimension $r$; 
\item there is an action of $\mu_N$ on $\mr{Hyp}$ for which $\pi$ is $\mu_N$-equivariant (acting trivially on the base);
\item after the base-change $\phi \colon S_{N, n, \ell} \to k$, the lisse sheaf
\[
\mc{F}= \mr{Gr}^W_r \left( \mr{R}^r \pi_! \Qlb \right)^{\chi}
\]
on $\mathbb{A}^1_{S_{N, n, \ell}}-\{T_1, \ldots, T_n\}$ becomes isomorphic to $\mc{F}_k$. Here $\mr{Gr}^W_{\bullet}$ denotes the weight filtration of \cite{deligne:weil2}.
\end{itemize}
\end{thm}
For example, if $k= \Qb$ and the points $\alpha_i$ lie in $\mathbb{A}^1(\Q(\zeta_N))$, this result descends the given rigid local system $\mc{F}_{\Qb}$ to an arithmetic local system on $U/\Q(\zeta_N)$.

To prove Theorem \ref{katzmotivic}, Katz introduces a `geometric' operation on sheaves, the middle convolution, that can be iterated (in combination with a simpler twisting operation) to generate all cohomologically rigid local systems starting from easily-understood rank one local systems. We now describe (one version of) the middle convolution. First, for any $r \geq 0$, let $\mathbb{A}(n, r+1)_{R_{N, \ell}}$ be the space
\[
\mathbb{A}(n, r+1)_{R_{N, \ell}}= \Spec R_{N, \ell}[T_1, \ldots, T_n, X_1, \ldots, X_{r+1}]\left[ \frac{1}{\prod_{i \neq j}(T_i-T_j) \prod_{a, j} (X_a- T_j) \prod_k (X_{k+1}-X_k)}\right].
\]
For $r=0$, we recover $\mathbb{A}(n, 1)_{R_{N, \ell}}= \mathbb{A}^1_{S_{N, n, \ell}}-\{T_1, \ldots, T_n\}$ (in this case the product over $k$ in the above definition is understood to be 1). To ease notation, we will suppress the $R_{N, \ell}$ subscript and simply write $\mathbb{A}(n, r+1)$. Now, let $\mr{Lisse}(N, n, \ell)$ be the category of lisse $\Qlb$-sheaves on $\mathbb{A}(n, 1)$. For any non-trivial character $\chi \colon \mu_N(R_{N, \ell}) \to \Qlb^\times$, incarnated as a character sheaf $\mc{L}_{\chi}$ on $(\mathbb{G}_{m})_{R_{N, \ell}}$, and then as the pull-back $\mc{L}_{\chi(X_2-X_1)}$ on $\mathbb{A}(n, 2)$ via the difference map
\[
d \colon \mathbb{A}(n, 2) \xrightarrow{X_2-X_1} (\mathbb{G}_{m})_{R_{N, \ell}},
\]
Katz defines two operations, na\"{i}ve and middle convolution by $\chi$,
\begin{align*}
\mr{NC}_{\chi} \colon \mr{Lisse}(N, n, \ell) \to \mr{Lisse}(N, n, \ell), \\
 \mr{MC}_{\chi} \colon \mr{Lisse}(N, n, \ell) \to \mr{Lisse}(N, n, \ell),
\end{align*}
via the following recipe. For all $r \geq 0$ and $1 \leq i \leq r+1$, there are projections 
\begin{align*}
\pr_i \colon \mathbb{A}(n, r+1) &\to \mathbb{A}(n, 1) \\
(T_1, \ldots, T_n, X_1, \ldots, X_{r+1}) &\mapsto (T_1, \ldots, T_n, X_i).
\end{align*}
The na\"{i}ve convolution of $\mc{F}$ by $\chi$ is given by ($r=1$ in the above notation)
\begin{equation}\label{naive}
\mr{NC}_{\chi}(\mc{F})= \mr{R}^1(\pr_{2, !}) \left( \pr_1^* \mc{F} \otimes d^* \mc{L}_{\chi} \right).
\end{equation}
From now on, we will abbreviate $\pr_1^* \mc{F} \otimes d^* \mc{L}_{\chi}$ by simply $\mc{F} \boxtimes \mc{L}_{\chi}$. Note that if $\mc{F}$ is mixed of weights $\leq w$, then $\mr{NC}_{\chi}(\mc{F})$ is mixed of weights $\leq w+1$. The middle convolution is then a variant that will take pure sheaves to pure sheaves. To define it, view $\pr_2$ as a relative $\mathbb{A}^1$ (with coordinate $X_1$) with the sections $T_1, \ldots, T_n, X_2$ deleted, and compactify it to a relative $\mathbb{P}^1$:
\[
\xymatrix{
\mathbb{A}(n, 2) \ar[d]^{\pr_1} \ar[dr]^{\pr_2} \ar[r]^-j & \mathbb{P}^1 \times_{R_{N, \ell}} \mathbb{A}(n, 1) \ar[d]^{\overline{\pr}_2} \\
\mathbb{A}(n, 1) & \mathbb{A}(n, 1).
}
\] 
For all $\mc{F} \in \mr{Lisse}(N, n, \ell)$, set
\begin{equation}\label{middle}
\mr{MC}_{\chi}(\mc{F})= \mr{R}^1 (\overline{\pr}_{2, *}) \left( j_*( \mc{F} \boxtimes \mc{L}_{\chi}) \right).
\end{equation}
See \cite[Lemma 8.3.2]{katz:rls} for the basic facts about $\mr{MC}_{\chi}(\mc{F})$: note for now that it is lisse, and for non-trivial $\chi$ and $\mc{F}$ pure of weight $w$, the natural map
\[
\mr{NC}_{\chi}(\mc{F}) \to \mr{MC}_{\chi}(\mc{F}) 
\]
is a surjection identifying $\mr{MC}_{\chi}(\mc{F})$ with the top graded piece $\Gr^W_{w+1}(\mr{NC}_{\chi}(\mc{F}))$ of the weight filtration on $\mr{NC}_{\chi}(\mc{F})$.

The proof of Theorem \ref{katzmotivic} is a combination of Katz's main theorem (see \cite[Theorem 5.2.1]{katz:rls}) on the structure of tamely ramified quasi-unipotent rigid local systems on $U/k$ ($k$ algebraically closed, recall), which shows that they are all obtained by starting from rank one local systems and iterating some combination of the middle convolution (specialized via $\phi \colon S_{N, n, \ell} \to k$) and the simpler operation of twisting by rank one local systems, and understanding the total geometric effect of iterating the convolution. 
Our main result will require invoking a variant of Katz's results. We will need to analyze local systems that are not necessarily rigid but that nevertheless admit a decription as an iteration of middle convolutions and twists (as in Katz), as well as Schur functors. We will need to check that the natural arithmetic descents of these local systems specialize (over number fields) to the $\ell$-adic realizations of motives for motivated cycles (in the sense of \cite{andre:motivated}). This will make precise the sense in which the Kuga-Satake lifts we construct are `motivic.' Our notation here is that for any object $M$ of $\mc{M}_{F, E}$, and any embedding $\iota \colon E \into \Qlb$, the corresponding $\ell$-adic realization (a $\gal{F}$-representation) is denoted $M_{\iota, \ell}$. We now explain all of this in detail.

Now we come to the key technical lemma, which is a simple variant of the arguments used by Katz in \cite[Chapter 8]{katz:rls} to establish Theorem \ref{katzmotivic}; Dettweiler has conveniently abstracted many of these arguments in \cite{dettweiler:middleconvolution}, and we will be able to reduce to these already-treated cases. To simplify the statement, we first make a definition:
\begin{defn}\label{concentrated}
An object $\mc{F}$ of $\mr{Lisse}(N, n, \ell)$ is said to be \textit{geometrically concentrated} if it is of the following form: there exists a smooth morphism $f \colon \mathbf{X} \to \mathbb{A}(n, 1)$ such that
\begin{itemize}
\item $f$ is equivariant for the action of a finite group $G$ on $\mathbf{X}$, trivial on the base $\mathbb{A}(n, 1)$; 
\item there is an idempotent $e \in \Qlb[G]$;
\item there is an integer $r \geq 0$ such that for all $i \neq r$, $e \mr{R}^if_! \Qlb=0$, and:
\item $\mc{F} \cong \Gr^W_r (e \mr{R}^r f_! \Qlb)$.
\end{itemize}
We will also use the terminology `geometrically concentrated' for sheaves on $\mathbb{A}^1-\{\alpha_1, \ldots, \alpha_n\}$ admitting a similar description, where we have specialized the parameters $T_i$ to some $\alpha_i \in \mathbb{A}^1(k)$.
\end{defn}
Here is the result:
\begin{prop}\label{geometricops}
Let $\mc{F}$ be any object of $\mr{Lisse}(N, n, \ell)$ obtained by some iteration, starting from the constant sheaf, of the following three operations:
\begin{itemize}
\item middle convolution by a non-trivial character $\rho \colon \mu_N(R_{N, \ell}) \to \Qlb^\times$;
\item Schur functors;
\item tensoring by geometrically concentrated sheaves in $\mr{Lisse}(N, n, \ell)$.
\end{itemize}
Then $\mc{F}$ is geometrically concentrated.

Moreover, let $F$ be a field of characteristic zero (eg, a number field) containing $R_{N, \ell}$, and let $t \colon \Spec F \to \mathbb{A}(n, 1)$ be any map of $R_{N, \ell}$-schemes. Then there exists a number field $E$ and, for each embedding $\iota \colon E \into \Qlb$, an object $M$ of $\mc{M}_{F, E}$ such that $t^* \mc{F}$, as $\gal{F}$-representation, is isomorphic to the $\ell$-adic realization $M_{\iota, \ell}$.
\end{prop}
\proof
Let $\mc{F}$ be geometrically concentrated. Arguing inductively, we have to check the following three statements:
\begin{itemize}
\item if $\rho \colon \mu_N(R_{N, \ell}) \to \Qlb^\times$ is a non-trivial character, then $\mr{MC}_{\rho}(\mc{F})$ is geometrically concentrated;
\item if $\lambda$ is a partition of a positive integer $m$, then the image of the associated Schur functor $\mathbb{S}_{\lambda} \mc{F}$ is geometrically concentrated;
\item and if $\mc{G}$ is a second geometrically concentrated object of $\mr{Lisse}(N, n, \ell)$, then $\mc{F} \otimes \mc{G}$ is geometrically concentrated.
\end{itemize}
The first assertion follows from \cite[2.6.1 Theorem]{dettweiler:middleconvolution}. For the second bulleted assertion, let $(\mc{F}, f, e, r)$ be as in Definition \ref{concentrated}, and let
\[
f^m \colon \mathbf{X} \times_{\mathbb{A}(n, 1)} \mathbf{X} \times_{\mathbb{A}(n, 1)} \cdots \times_{\mathbb{A}(n, 1)} \mathbf{X} \to \mathbb{A}(n, 1)
\]
denote the projection from the $m$-fold fiber product. This projection is now equivariant for the natural action of $G^m \rtimes S_m$. By the K\"{u}nneth formula and the assumption of concentration,
\[
(e \times e \times \cdots \times e) \mr{R}^{mr} f^m_! \Qlb\cong (e \times e \times \cdots \times e) \bigoplus_{i_1 + \cdots +i_m= mr} \left( \otimes_{j=1}^m \mr{R}^{i_j}f_! \Qlb\right) \cong (e \mr{R}^r f_! \Qlb)^{\otimes m}.
\]
Applying $\Gr^W_{mr}$, and abbreviating $e^m= e \times \cdots \times e$, we deduce that
\[
e^m \Gr^W_{mr} \mr{R}^{mr} f^m_! \Qlb \cong \mc{F}^{\otimes m}.
\]
Let $s_{\lambda} \in \Qlb[S_m]$ be the idempotent projector such that for any vector space $V$, $\mathbb{S}_{\lambda}(V)= s_{\lambda}(V^{\otimes m})$. Then the element $s_{\lambda}\cdot e^m$ of $\Qlb[G^m \rtimes S_m]$ is still idempotent, since the operators $s_{\lambda}$ and $e^m$ commute. Thus
\[
\mathbb{S}_{\lambda}(\mc{F})= s_{\lambda} e^m \left( \Gr^W_{mr} (\mr{R}^{mr}f^m_! \Qlb) \right)
\]
is geometrically concentrated.

For the last bulleted point, again apply the K\"{u}nneth formula, using concentration of $\mc{F}$ and $\mc{G}$.

The second half of the proposition-- the claim that for all specializations $t$, $t^* \mc{F}$ is the $\ell$-adic realization of a motivated motive-- is deduced from the first part using the fact that for any smooth (not necessarily projective) variety $U/F$, and for all $k \geq 0$, the top graded quotient $\Gr^W_k H^k_c(U_{\overline{F}}, \Qlb)$ can be conveniently described in terms of a smooth compactification of $U$ (with a smooth normal crossings divisor at the boundary). For details, see the discussion in between Remark 2.6 and Corollary 2.7 in \cite{patrikis-taylor:irr}.
\endproof
\begin{rmk}
The following rank 1 geometrically concentrated sheaves are, along with the middle convolution, the ingredients in Katz's `motivic' construction of rigid local systems: fix $n$ characters
\[
\chi_{i} \colon \mu_N(R_{N, \ell}) \to \Qlb^\times,
\]
$i=1, \ldots, n$, $a=1, \ldots, r+1$, and set 
\[
\mc{F}= \otimes_{i=1}^n (X-T_i)^*(\mc{L}_{\chi_{i}}),
\]
where $X$ denotes the `relative $\mathbb{A}^1$' coordinate on $\mathbb{A}(n, 1)$, whence a map $(X-T_i) \colon \mathbb{A}(n, 1) \to \mathbb{G}_m$. All of the `universal' extensions of rigid local systems to $\mr{Lisse}(N, n, \ell)$ are then expressed by iterating middle convolution and tensoring by such rank 1 sheaves $\mc{F}$.
\end{rmk}
\section{The main result}
We will now demonstrate the existence of `generalized Kuga-Satake lifts' of motives whose $\ell$-adic realizations are isomorphic to some fiber of an $\mr{SO}_5$ or $\mr{SO}_6$-cohomologically rigid local system. Note that an $\mr{SO}_N$-local system that is $\mr{GL}_N$-cohomologically rigid is automatically $\mr{SO}_N$-cohomologically rigid: $H^1(\mathbb{P}^1, j_* \mf{so}_N)$ is a summand of $H^1(\mathbb{P}^1, j_* \mf{gl}_N)$, so vanishing of the latter group of course implies vanishing of the former group. In light of Scholium \ref{rigimpliesmot}, the essential content of our result is the combination of Katz's theory with the following theorem of Bogner-Reiter; what we will do is reinterpret this theorem and give a couple examples:
\begin{thm}[Theorem 3.1 of \cite{bogner-reiter:rigidsp4}]\label{bognerreiter}
Each $\mr{Sp}_4(\CC)$-cohomologically rigid local system with quasi-unipotent local monodromies comes from geometry: it can be constructed by a sequence of `geometric operations' (direct sum, tensor product, exterior square, symmetric square, rational pull-back, and middle convolution), starting from rank 1 local systems.
\end{thm}
The argument crucially depends on `low-dimensional coincidences' in the Dynkin diagrams of the classical Lie groups. To give a concrete sense of the difficulties involved in generalizing to higher rank, we will construct (in Examples \ref{so7} and \ref{so7bis} below) motivic, $\mr{SO}_7$-rigid local systems whose $\mr{Spin}_7$ lifts do not seem to be provably motivic using the methods of this paper.
\begin{thm}\label{main}
Let $F$ and $E$ be number fields, and fix an embedding $\iota \colon E \into \Qlb$. Suppose $M$ is a motivated motive in $\mc{M}_{F, E}$ whose $(\iota, \ell)$-adic realization is orthogonal of rank 5, and is isomorphic to a $\gal{F}$-representation of the following form: there exists 
\begin{itemize}
\item a finite set of points $S \subset \mathbb{P}^1(F)$, with complement $U= \mathbb{P}^1-S$;
\item an orthogonal rank 5 local system $\mc{F}$ on $U$ over $F$, such that $\mc{F}_{\overline{F}}$ (the corresponding geometric local system) has quasi-unipotent local monodromies and is $\mr{SO}_5$-cohomologically rigid;
\item and a point $t \colon \Spec F \to U$;
\item such that $M_{\iota, \ell}$ is isomorphic to $t^* \mc{F}$ as $\Qlb[\gal{F}]$-module.
\end{itemize}
For simplicity, assume either that $M_{\iota, \ell}$ is irreducible, or that the geometric monodromy group of $\mc{F}_{\overline{F}}$ acts irreducibly in the standard 5-dimensional representation. Then (enlarging $F$ if necessary by a quadratic extension so that $M_{\iota, \ell}$ is special orthogonal) there exists a geometric lift
\[
\xymatrix{
& \mr{GSpin}_5(M_{\iota, \ell}) \ar[d] \\
\gal{F} \ar@{-->}^{\tilde{\rho}_t}[ur] \ar[r] & \mr{SO}_5(M_{\iota, \ell})
}
\]
such that for some finite extensions $F'/F$ and $E'/E$ of number fields, $r_{spin} \circ \tilde{\rho}|_{\gal{F'}}$ is isomorphic to the $E' \xrightarrow{\iota'} \Qlb$-realization of an object of $\mc{M}_{F', E'}$, for a suitable embedding $\iota'$ extending $\iota$.
\end{thm}
\begin{rmk}
The hypothesis that either $M_{\iota, \ell}$ or $\mc{F}_{\overline{F}}$ has `big' monodromy is to avoid (more) tedious case-by-case arguments; it will be met in the examples discussed below.
\end{rmk}
\proof
The first part of the proof is just bookkeeping to reduce to the case where $\mc{F}_{\overline{F}}$ has monodromy group containing $\mr{SO}_5$; the reader may wish simply to skip to that part of the proof (fourth paragraph). If $M_{\iota, \ell}$ is irreducible, then by \cite[Proposition 3.4.1]{stp:variationsarxiv} it can be written in the form $\Ind_{\gal{L}}^{\gal{F}}(r \otimes \tau)$ where $r$ is Lie-irreducible (irreducible after every finite base-change) and $\tau$ has finite image. Since the rank of $M_{\iota, \ell}$ is 5, this means $M_{\iota, \ell}$ is either the induction of a character from a rank 5 extension $L/F$, or is Lie-irreducible, or has finite-image. The finite-image case is obvious, since then by Tate's theorem (\cite[Theorem 4]{serre:DSsurvey}) there exists a $\tilde{\rho}$ with finite image, for which $r_{spin} \circ \tilde{\rho}$ is the realization of an Artin motive. The induction case follows since any geometric Galois character $\gal{L} \to \Qlb^\times$, at least after a finite base-change $L'/L$ can be cut out as a motivated sub-motive of a CM abelian variety. It is easy to see that the same will then (after a finite base-change) be true of any geometric lift $\tilde{\rho}$ (and such geometric lifts are easy to produce in this case). So we may assume $M_{\iota, \ell}$ is Lie-irreducible. There are then only two possibilities for the connected component of the Zariski closure of the image of this $\gal{F}$-representation: either it is the image of the symmetric fourth power of $\mr{SL}_2$, or it is all of $\mr{SO}_5$. 

Let us denote by $\rho \colon \pi_1(U, \bar{t}) \to O(\mc{F}_{\bar{t}}) \cong O_5(\Qlb)$ the representation associated to the lisse sheaf $\mc{F}$ (with $\bar{t}$ the geometric point over $t$ corresponding to our choice of algebraic closure $\overline{F}/F$). If in the last paragraph, $M_{\iota, \ell}$ had (after passing to the connected component) $\mr{SO}_5$ as its monodromy group, then of course $\overline{\rho(\pi_1(U, \bar{t}))}^{Zar}$ contains $\mr{SO}_5$. In particular, the geometric monodromy group $\overline{\rho(\pi_1(U_{\overline{F}}, \bar{t}))}^{Zar}$ has connected component equal to a connected normal subgroup of $\mr{SO}_5$, i.e. either $\mr{SO}_5$ itself or the trivial group. If trivial, then the finite group $\overline{\rho(\pi_1(U_{\overline{F}}, \bar{t}))}$ is normalized by the connected group $\overline{\rho(\pi_1(U, \bar{t}))}^{Zar,0}= \mr{SO}_5$; but of course any map (conjugation, here) from a connected group to a finite group is trivial, so $\mr{SO}_5$ centralizes $\overline{\rho(\pi_1(U_{\overline{F}}, \bar{t}))}$, i.e. the latter group is contained in $Z(O_5)= \{\pm 1\}$. But such a local system cannot be rigid.

Next, if $M_{\iota, \ell}$ has monodromy group $\mr{PGL}_2$ acting through $\Sym^4$, then the connected component $\overline{\rho(\pi_1(U_{\overline{F}}, \bar{t}))}^{Zar, 0}$ is either $\{1\}$ or $\mr{PGL}_2$ or $\mr{SO}_5$. By the argument of the previous paragraph, it cannot be trivial, so let us consider the case of $\mr{PGL}_2$. The normalizer of $\mr{PGL}_2$ in $\mr{SO}_5$ is simply $\mr{PGL}_2$,\footnote{More generally, one can identify the normalizer in a simple algebraic group of the principal $\mr{SL}_2$, i.e. the output of the Jacobson-Morosov theorem for a regular nilpotent element. For example, for $G= \mr{SO}_{2n+1}$, it is easy to see that the centralizer $Z_G(\mr{PGL}_2)=\{1\}$, and since $\mr{PGL}_2$ has no outer automorphisms, we deduce that $N_G(\mr{PGL}_2)= \mr{PGL}_2$.} so $\overline{\rho(\pi_1(U, \bar{t}))}^{Zar}$ must equal either $\mr{PGL}_2$ or $\mr{PGL}_2 \times \{\pm 1\}$. In this case a considerably simpler version of the argument we give below will apply, solving an $\mr{SL}_2 \to \mr{PGL}_2$ lifting problem rather than a $\mr{Spin}_5 \to \mr{SO}_5$ lifting problem; we omit the details, noting only that when we lift the geometric local system $\mc{F}_{\overline{F}}$ to $\mr{SL}_2$, we obtain a linearly rigid rank 2 local system, which can be described by Katz's algorithm.

Finally, we treat the basic case-- this is the heart of the theorem-- in which the geometric monodromy group of $\mc{F}_{\overline{F}}$ is irreducible. We make one initial adjustment to $\rho$-- denoting as before the monodromy representation of $\mc{F}$-- which will eventually cost us at most a quadratic extension on the field $F'$ in the conclusion of the theorem. Namely, $\rho$ may have non-trivial determinant 
\[
\det \rho \colon \pi_1(U, \bar{t}) \to \{ \pm 1 \}.
\]
The underlying geometric (rank 1) local system $\det \rho|_{\pi_1(U_{\overline{F}}, \bar{t})}$, if non-trivial, is isomorphic to a tensor product of translated Kummer sheaves (see \cite[\S 8.1]{katz:rls} for this notion) of order 2. Such translated Kummer sheaves descend to $\pi_1(U, \bar{t})$, so there exists a lisse character sheaf $\mc{L}$ on $U$ over $F$ of order at most 2 such that $\det( \rho \otimes \mc{L})$ is geometrically trivial, i.e. factors through a character $\gal{F} \to \{ \pm 1\}$. Thus, there exists an at most quadratic extension $F'/F$ such that $\det(\rho \otimes \mc{L}) =1$ as $\pi_1(U_{F'}, \bar{t})$-representation. From now on, we replace $F$ by $F'$ and $\rho$ (respectively, $\mc{F}$) by $\rho \otimes \mc{L}$ (respectively, $\mc{F} \otimes \mc{L}$), so we can work in the setting of the following diagram:
\begin{equation}\label{twisttoso}
\xymatrix{
& \mr{Spin}(\mc{F}_{\bar{t}}) \ar[r]^{r_{spin}}_{\sim} \ar[d]^{\pi} & \mr{Sp}(W_{spin, \bar{t}}) \ar[r] & \mr{GL}(W_{spin, \bar{t}}) \\
\pi_1(U, \bar{t}) \ar[r]_{\rho} \ar@{-->}[ur]^{???} & \mr{SO}(\mc{F}_{\bar{t}}), & {} \\
}
\end{equation}
where we let $W_{spin, \bar{t}}$ denote the spin representation, with its natural (in this case symplectic) bilinear pairing-- this pairing is well-defined up to $\Qlb^\times$-scaling. For the time being we set aside the arithmetic local system $\mc{F}$ and only work with $\mc{F}_{\overline{F}}$. 

The representation $\rho$ lands in a profinite subgroup of $\mr{SO}(\mc{F}_{\bar{t}})$ (namely, a group isomorphic to $\mr{SO}_5(\mc{O})$ for the ring of integers $\mc{O}$ in some finite extension of $\Ql$). We claim that the restriction $\rho|_{\pi_1(U_{\overline{F}}, \bar{t})}$ lifts to $\mr{Spin}(\mc{F}_{\bar{t}})$; but we can obviously lift the restriction of $\rho$ to the free group (not yet profinitely-completed) on $|S|-1$ elements to a homomorphism landing in a profinite subgroup of $\mr{Spin}(\mc{F}_{\bar{t}})$, and so the universal property of profinite completion implies the claim. Let us call such a lift $\tilde{\rho}_{\overline{F}}$.

The Theorem \ref{bognerreiter} of Bogner and Reiter now says that there is a local system $\mc{G}_{\overline{F}}$ on $U_{\overline{F}}$, constructed by a series of `geometric' operations, that is isomorphic as rank 4 local system to $r_{spin} \circ \tilde{\rho}_{\overline{F}}$. Choose an isomorphism, so that we can identify $\rho_{\mc{G}_{\overline{F}}}$ and $r_{spin} \circ \tilde{\rho}_{\overline{F}}$ as homomorphisms to $\mr{GL}(W_{spin, \bar{t}} \xrightarrow{\sim} \mc{G}_{\overline{F}, \bar{t}})$. Note that the necessary `geometric operations' are, starting from rank one sheaves (all of which are products of translated Kummer sheaves): middle convolution by non-trivial characters, application of Schur functors, tensor product of two sheaves obtained through the previous operations-- note that these all fall under the umbrella of Proposition \ref{geometricops}-- and finally pull-backs along rational maps $\mathbb{P}^1 \to \mathbb{P}^1$. We have not taken these rational pull-backs into account in Proposition \ref{geometricops}, and we do so now in an \textit{ad hoc} manner, working one-by-one through the relevant cases in the proof of \cite[Theorem 3.1]{bogner-reiter:rigidsp4}.\footnote{The reason anything more needs to be said at this point is that pullback is not, strictly speaking, a `reversible operation': applying a series of middle convolutions, twists, and pull-backs to our given $r_{spin} \circ \tilde{\rho}_{\overline{F}}$ to produce something visibly geometric does not immediately allow us to express $r_{spin} \circ \tilde{\rho}_{\overline{F}}$ as a `geometrically concentrated' (Definition \ref{concentrated}) sheaf.} We will use freely the notation of that argument; the only cases then needing to be considered are $P_3(4, 8, 10, 10)$, $P_4(6, 6, 10, 10)$, $P_4(6, 8, 8, 10)$, and $P_5(8, 8, 8, 10, 10)$ (the other cases of their theorem making use of pullbacks all correspond to \textit{reducible} $\mr{SO}_5$-local systems, which we have by assumption excluded).
\begin{itemize}
\item The cases $P_3(4, 8, 10, 10)$, $P_4(6, 6, 10, 10)$, $P_4(6, 8, 8, 10)$ are all dealt with similarly; for these, \cite[Theorem 3.1]{bogner-reiter:rigidsp4} shows that a series of middle convolution and tensoring by rank 1 sheaves reduces our given local system $r_{spin} \circ \tilde{\rho}_{\overline{F}}$ to a $\mr{GO}_2$-local system. Any such sheaf $\mc{H}$ is the induction from an index two subgroup of a rank 1 local system $\mc{L}$, i.e. there exists a finite \'{e}tale (degree 2) cover $\pi \colon U' \to U$ such that $\mc{H} \cong \pi_* \mc{L}$. This is of course geometrically concentrated, so applying (in reverse) the relevant middle convolutions and tensor products by Kummer sheaves we obtain a geometrically concentrated construction of $r_{spin} \circ \tilde{\rho}_{\overline{F}}$.
\item In the remaining case, $P_5(8, 8, 8, 10, 10)$, \cite[Theorem 3.1]{bogner-reiter:rigidsp4} yields an identity of the form 
\[
\mc{F}_1 \otimes \mc{F}_2 \cong f^* \left(\mr{MC}_{-1}(r_{spin} \circ \tilde{\rho}_{\overline{F}}) \right),
\] 
where $f \colon \mathbb{P}^1 \to \mathbb{P}^1$ is a degree 2 map, and the $\mc{F}_i$ are $\mr{GL}_2$-rigid. Restricting $f$ to an \'{e}tale, necessarily Galois, cover $U' \to U$, and applying Frobenius reciprocity, we find that 
\[
f_*(\mc{F}_1 \otimes \mc{F}_2) \cong \mr{MC}_{-1}(r_{spin} \circ \tilde{\rho}_{\overline{F}}) \oplus \left( \mr{MC}_{-1}(r_{spin} \circ \tilde{\rho}_{\overline{F}}) \otimes \delta_{U'/U} \right),
\]
where $\delta_{U'/U}$ denotes the non-trivial quadratic character of $\pi_1(U_{\overline{F}})/\pi_1(U'_{\overline{F}})$. We can then construct an idempotent $e \in \Qlb[ \Aut(U'/U)]$ such that
\[
e \left( f_*(\mc{F}_1 \otimes \mc{F}_2) \right) \cong \mr{MC}_{-1}(r_{spin} \circ \tilde{\rho}_{\overline{F}}).
\] 
It follows that $r_{spin} \circ \tilde{\rho}_{\overline{F}}$ admits a description as a geometrically concentrated local system.
\end{itemize}

The explicit construction of $\mc{G}_{\overline{F}}$ gives us more than just a geometric local system, however. The discussion of \S \ref{rigid} shows that for some cyclotomic extension $K$ of $F$, depending only on the geometric local monodromies of $\rho$,\footnote{If all eigenvalues of local monodromy of $\rho$ are $N^{th}$ roots of unity, we can take $F(\mu_{2N})$.} there exists an arithmetic descent $\mc{G}$ of $\mc{G}_{\overline{F}}$ to a lisse sheaf on $U_K$. We denote by $\tilde{\rho}$ the corresponding representation of $\pi_1(U_K, \bar{t})$. It takes values in the normalizer of the geometric monodromy group of $\mc{G}_{\overline{F}}$, which is necessarily contained in $\mr{GSp}(W_{spin, \bar{t}})$. Since we can extend $\pi$ (technically, $\pi \circ r_{spin}^{-1}$) to a quotient map (killing the center) $\pi \colon \mr{GSp}(W_{spin, \bar{t}}) \to \mr{SO}(W_{spin, \bar{t}})$, we can now compare $\pi \circ \tilde{\rho}$ with $\rho$ as $\pi_1(U_K, \bar{t})$-representations. The space
\[
\Hom_{\pi_1(U_{\overline{F}}, \bar{t})}(\pi \circ \tilde{\rho}, \rho)
\]
is a one-dimensional $\Qlb$-vector space with an action of $\pi_1(U_K, \bar{t})/\pi_1(U_{\overline{K}}, \bar{t}) \cong \gal{K}$. The fact that both $\pi \circ \tilde{\rho}$ and $\rho$ land in $\mr{SO}_5$ forces this action to be trivial.\footnote{Letting $\chi$ denote the resulting character of $\gal{K}$, $\chi^2$ and $\chi^5$ both equal 1, by comparing orthogonal multipliers and determinants.} Thus in fact $\pi \circ \tilde{\rho} = \rho$ as $\pi_1(U_K, \bar{t})$-representation (they are isomorphic, and we already knew they were equal on $\pi_1(U_{\overline{F}}, \bar{t})$). In particular, pulling back along $t \colon \Spec K \to U_K$ (the restriction of the original $t \in U(F)$), we have a commutative diagram of $\gal{K}$-representations
\[
\xymatrix{
& \mr{GSpin}(\mc{F}_{\bar{t}}) \ar[d]^\pi \\
\gal{K} \ar[ur]^{\tilde{\rho}_t} \ar[r]_-{\rho_t} & \mr{SO}(\mc{F}_{\bar{t}}).
}
\]
Crucially, by Proposition \ref{geometricops}, $\tilde{\rho}_t$ (which by definition is $\rho_{\mc{G}, t}$) is the $\ell$-adic realization of a motivated motive. We have therefore produced a generalized Kuga-Satake lift of the original motive $M$, after some base-change: the $\ell$-adic realization (now restricted to $\gal{K}$) of $M$ was assumed to be isomorphic to $\rho_t$.
\endproof
\begin{rmk}\label{inefficiency}
Let us record how inefficient the argument was in preserving the original field of definition $F$ of the motive $M$, in the most important case of the theorem in which the geometric monodromy group of $\mc{F}_{\overline{F}}$ was irreducible. We were given a local system $\mc{F}$ with eigenvalues of local monodromy generating some cyclotomic field $\Q(\mu_N)$; after replacing $\mc{F}$ by a quadratic twist $\mc{F} \otimes \mc{L}$, we (maybe) had to enlarge $F$ by a quadratic extension to make $\mc{F}$ a special orthogonal local system; the lifted arithmetic local system $\widetilde{\mc{F}}$ was defined over $F(\mu_{2N})$; and finally, taking into account the original twist $\mc{F} \mapsto \mc{F} \otimes \mc{L}$, we had to replace $F$ by possibly one further quadratic extension. In sum, we have passed from $F$ to $F'(\mu_{2N})$, where $F'$ is some (at most) biquadratic extension, and $N$ is determined by the local monodromies of the local system $\mc{F}$. 
\end{rmk}
We will now make this procedure explicit for motivic lifting of fibers of the sextic Dwork family, establishing Theorem \ref{dwork}:
\begin{eg}
Recall that we are considering motives of the form $M_t= \Prim^4(\mathbf{X}_t)^\Gamma(2)$, where 
\[
\mathbf{X} \xrightarrow{f} \mathbb{P}^1-\{\infty, \mu_6\}
\]
is the smooth projective family with $t$-fiber
\[
\mathbf{X}_t= \{\sum_{i=0}^5 X_i^6= 6t \prod_{i=0}^5 X_i \} \into \mathbb{P}^5,
\]
acted on by the finite group scheme $\Gamma$ given by the kernel of the multiplication map $\mu_6^6 \xrightarrow{mult} \mu_6$. $M_t$ is a well-defined object of $\mc{M}_{\Q(t)}$. For any $\ell$, the lisse $\Ql$-sheaf of $\Gamma$-invariants in primitive cohomology, 
\[
\mc{F}= (\mr{R}^4f_* \Ql)^{\Gamma}_{prim},
\] 
on $\mathbb{P}^1-\{\infty, \mu_6\}$ over $\Q$ has $t^* \mc{F} \cong M_{t, \ell}$ as $\gal{\Q(t)}$-representation. The geometric local system $\mc{F}_{\Qb}$ is not rigid, but it is closely related to a rigid local system on $\mathbb{P}^1-\{0, 1, \infty\}$. Namely, there exists a smooth projective family $\mathbf{X}' \xrightarrow{f'} \mathbb{P}^1-\{0, 1, \infty\}$ with $\Gamma$-action, and a $\Gamma$-equivariant pullback diagram
\[
\xymatrix{
\mathbf{X} \ar[d]^f \ar[r] & \mathbf{X}' \ar[d]^{f'} \\
\mathbb{P}^1-\{\infty, 0, \mu_6\} \ar[r]^{[6]} & \mathbb{P}^1-\{0, 1, \infty\},
}
\]
with $[6]$ denoting the map $z \mapsto z^6$.\footnote{There are various ways of producing such an $\mathbf{X}'$, for instance taking $\mathbf{X}'$ to be the family, with parameter $s$, 
\[
s^{-1}Y_0^6 +\sum_{i=1}^5 Y_i^6= 6 \prod_{i=0}^5 Y_i.
\]
This is observed by Katz in \cite{katz:dwork}.
}
Now the lisse $\Qlb$-sheaf $\mc{F'}= (\mr{R}^4 f'_* \Qlb)^\Gamma_{prim}$ satisfies $[6]^* \mc{F}' \cong \mc{F}$, hence setting $M'_t= \Prim^4(X'_t)^\Gamma(2)$, we have
\[
M_{t, \ell} \cong M'_{t^6, \ell},
\]
so for our purposes we can work with the motives $M'_t$ instead. Now, $\mc{F}'_{\Qb}$ is a $\mr{GL}_5$-cohomologically (hence physically) rigid local system satisfying an orthogonal autoduality (Poincar\'{e} duality). Its local monodromies are regular unipotent (at $\infty$), a reflection (at 1), and regular semi-simple (at 0) with eigenvalues $\{\zeta_6^i\}_{i=1, \ldots, 5}$. As in the proof of Theorem \ref{main}, we twist $\mc{F}'$ by the product of translated Kummer sheaves on $\mathbb{P}^1-\{0,1,\infty\}$ over $\Q$ having geometric local monodromies $-1$ at 0 and 1. The result is an $\mr{SO}_5$-local system (at least geometrically) whose lift $\widetilde{\mc{F}}$ to $\mr{Sp}_4$ has local monodromies with Jordan forms
\[
U(4) \qquad -1,-1,1,1 \qquad \eta, \eta^3, \eta^9, \eta^{11},
\]
where $U(4)$ indicates a 4-by-4 unipotent Jordan block and $\eta^2= \zeta_6$. Note that the sum of the centralizer dimensions is $4+8+4= 16$, so this local system is not $\mr{GL}_4$-cohomologically rigid. The Bogner-Reiter prescription in this case allows us to produce a rank 4 local system with these monodromies  as follows. We first compute $\mr{MC}_{-1}(\widetilde{\mc{F}})$, where $-1$ denotes the non-trivial character $\mu_2 \to \Qlb^\times$. This has local monodromies, in the above notation,
\[
-U(5), -1 \qquad -\eta, -\eta^3, -\eta^9, -\eta^{11}, 1, 1 \qquad U(2), U(2), 1, 1.
\]
After twisting (this is not essential) by the rank 1 sheaf with local monodromies $-1$ at 0 and $\infty$, we note that this $\mr{SO}_6$ local system lifts (via $\wedge^2$) to an $\mr{SL}_4$-local system with local monodromies
\[
U(4) \qquad \eta^2, \eta^4, \eta^7, \eta^{11} \qquad U(2), 1, 1.
\]
The sum of the centralizer dimensions here is $4+4+10=18$, so this rank 4 local system is cohomologically rigid, and therefore can be constructed (geometrically) by Katz's algorithm. In summary, we can run this procedure in reverse to construct, via the allowable geometric operations, a rank 4 symplectic local system $\mc{G}$ isomorphic to $\widetilde{\mc{F}}$. Thus we have made explicit, in the case of the sextic Dwork family, the procedure of Theorem \ref{main}; for a precise geometric description of these Kuga-Satake lifts, we need only combine Proposition \ref{geometricops} with \cite[Theorem 8.4.1]{katz:rls}, where the variety denoted $\mr{Hyp}$ in Theorem \ref{katzmotivic} is described as an explicit affine hypersurface.
\end{eg}
Let us also remark, without lingering over any details, that we similarly (but more easily) obtain a Kuga-Satake lifting result for motives whose $\ell$-adic realizations arise as fibers of certain $\mr{SO}_6$-local systems:
\begin{thm}\label{so6}
Let $F$ and $E$ be number fields, and fix an embedding $\iota \colon E \into \Qlb$. Suppose $M$ is a motivated motive in $\mc{M}_{F, E}$ whose $(\iota, \ell)$-adic realization is orthogonal of rank 6, and is isomorphic to a $\gal{F}$-representation of the following form: there exists 
\begin{itemize}
\item a finite set of points $S \subset \mathbb{P}^1(F)$, with complement $U= \mathbb{P}^1-S$;
\item a local system $\mc{F}$ on $U$ over $F$ such that $\mc{F}_{\overline{F}}$ has quasi-unipotent local monodromies and is $\mr{SO}_6$-cohomologically rigid (for example, $\mr{GL}_6$-cohomologically rigid) and irreducible;
\item and a point $t \colon \Spec F \to U$;
\item such that $M_{\iota, \ell}$ is isomorphic to $t^* \mc{F}$ as $\Qlb[\gal{F}]$-module.
\end{itemize}
Then (enlarging $F$ if necessary by a quadratic extension so that $M_{\iota, \ell}$ is special orthogonal) there exists a geometric lift
\[
\xymatrix{
& \mr{GSpin}_6(M_{\iota, \ell}) \ar[d] \\
\gal{F} \ar@{-->}^{\tilde{r}}[ur] \ar[r] & \mr{SO}_6(M_{\iota, \ell})
}
\]
such that for some finite extensions $F'/F$ and $E'/E$ of number fields, $r_{spin} \circ \tilde{r}|_{\gal{F'}}$ is isomorphic to the $E' \xrightarrow{\iota'} \Qlb$-realization of an object of $\mc{M}_{F', E'}$, for a suitable embedding $\iota'$ extending $\iota$. (Recall that $r_{spin}$ is now the sum of half-spin representations.)

\end{thm}
\proof
(Sketch) Let $\rho \colon \pi_1(U, \bar{t}) \to O_6(\Qlb)$ denote the monodromy representation of $\mc{F}$, with restriction $\rho_{\overline{F}}$ to $\pi_1(U_{\overline{F}}, \bar{t})$. The relevant low-dimensional coincidence is $\mr{Spin}_6 \cong \mr{SL}_4$, so we begin by choosing a lift
\[
\xymatrix{
& \mr{SL}_4(\Qlb) \ar[d]^{\pi} \\
\pi_1(U_{\overline{F}}, \bar{t}) \ar@{-->}[ur]^{\tilde{\rho}_{\overline{F}}} \ar[r]_{\rho_{\overline{F}}} & \mr{SO}_6(\Qlb).
}
\]
Since $\rho_{\overline{F}}$ is $\mr{SO}_6$-cohomologically rigid, $\tilde{\rho}_{\overline{F}}$ is necessarily (by equation \ref{cohrig}) $\mr{SL}_4$-cohomologically rigid, hence $\mr{GL}_4$-cohomologically rigid. Theorem \ref{katzmotivic} yields a geometric construction of $\tilde{\rho}$ via iterated middle convolution and tensoring with translated Kummer sheaves. This also yields an arithmetic descent $\tilde{\rho}$ of $\tilde{\rho}_{\overline{F}}$, and by irreducibility of $\rho_{\overline{F}}$, we can make the same argument as in Theorem \ref{main} to show that $\pi \circ \tilde{\rho}= \rho$ as $\pi_1(U_{F'}, \bar{t})$-representation, for some finite extension $F'/F$. We conclude the proof as in Theorem \ref{main}.
\endproof
We conclude with some examples indicating the difficulty of generalizing the arguments of this note:
\begin{eg}\label{so7}
With the Dwork family at hand, it makes sense to ask, for any even $n$, whether we can understand Kuga-Satake lifts of the orthogonal motives 
\[
\Prim^{n-2}(X_t)^{\Gamma}(\frac{n-2}{2})
\]
arising from the family of varieties $\mathbf{X} \xrightarrow{f} \mathbb{P}^1-\{\infty, \mu_n\}$ with $t$-fiber
\[
X_0^n+X_1^n+ \cdots+ X_{n-1}^n= ntX_0\cdots X_{n-1},
\]
carrying as before the natural action of $\Gamma= \ker(\mu_n^n \xrightarrow{mult} \mu_n)$. As for $n=6$, there is a $\mr{GL}_{n-1}$-cohomologically rigid local system, on $\mathbb{P}^1-\{0, 1, \infty\}$ whose pullback under $[n]$ is isomorphic to the lisse sheaf $(\mr{R}^{n-2}f_* \Qlb)^{\Gamma}_{prim}$. As before, we can twist to local system $\rho \colon \pi_1(\mathbb{P}^1-\{0,1,\infty\}) \to \mr{SO}_{n-1}(\Qlb)$, and ask: after lifting $\rho$ to $\tilde{\rho} \colon \pi_1(\mathbb{P}^1-\{0, 1, \infty\}) \to \mr{Spin}_{n-1}(\Qlb)$, can we find a `geometric' construction of a local system with local monodromies equal to those of $r_{spin} \circ \tilde{\rho}$? I do not know how to do this for $n \geq 8$. An indicator of the difficulty, even for $n=8$, using the methods of this paper might be the following: the spin representation factors $r_{spin} \colon \mr{Spin}_{7} \to \mr{SO}_8 \subset \mr{GL}_8$. The argument for $n=6$ crucially used the fact that composition with $r_{spin} \colon \mr{Spin}_5 \to \mr{Sp}_4$ yields a symplectic local system $r_{spin} \circ \tilde{\rho}$ that is (tautologically) $\mr{Sp}_4$-cohomologically rigid; so we might ask whether for $n=8$ $r_{spin} \circ \tilde{\rho}$ is $\mr{SO}_8$-cohomologically rigid? (\textit{A priori} it is only $\mr{Spin}_7$-rigid.) The answer is no: the Jordan forms of the local monodromies are
\[
U(7), 1 \qquad \imath^{\oplus 4}, -\imath^{\oplus 4} \qquad \zeta_8^{\pm 1}, \zeta_8^{\pm 2}, \zeta_8^{\pm 3}, 1, 1,
\]
for which the sum of the centralizer dimensions (in $\mr{SO}_8$) is $4+16+4=24 < 28= \dim \mf{so}_8$. Thus we will need in some way to come to terms with $\mr{Spin}_7$-rigidity.
\end{eg}
Is there any hope at all for tackling such problems? At least in some limited cases, I believe the answer is yes. I will close with one more $\mr{SO}_7$ example, which is $\mr{SO}_7$-cohomologically but not $\mr{GL}_7$-cohomologically rigid, and which admits a geometric construction. I expect as a result of work in progress to be able to show that the $\mr{Spin}_7$-lifts of this local system are (in the spin representation) geometric; indeed for every $n$ there will be an analogous example for $\mr{Spin}_{2n+1} \to \mr{SO}_{2n+1}$ lifts, so this should eventually provide examples of generalized Kuga-Satake theory in arbitrarily large rank.\footnote{For examples of a similar flavor, see \cite{stp:KSRLS2:heckeeigen}.}
\begin{eg}\label{so7bis}
There exists an $\mr{SO}_7$-cohomologically rigid local system on $U= \mathbb{P}^1-\{0, 1, \infty\}$ with local monodromies (at, in order, $\infty, 0, 1$):
\[
U(7) \qquad U(3), U(2), U(2) \qquad 1, -U(2), -U(2), -1, -1.
\]
Here is a construction. For a pair of characters $\alpha, \beta \colon \mu_N(R_{N, \ell}) \to \Qlb^\times$, with associated Kummer sheaves $\mc{L}_{\alpha}$ and $\mc{L}_{\beta}$, let us introduce the short-hand $\mc{L}(\alpha, \beta)$ for the local system on $U$ given by $\mc{L}_{\alpha} \otimes (x-1)^* \mc{L}_{\beta}$ (product of translated Kummer sheaves, so that $\alpha$ is placed at 0 and $\beta$ is placed at 1). Let $\pi$ denote the quotient map $\pi \colon \mr{Sp}_4 \to \mr{SO}_5$ (for an implicit choice of pairing that will arise below). Then the sheaf
\[
\mc{L}(1, -1) \otimes \mr{MC}_{-1} \left( \mc{L}(-1, 1) \otimes \mr{MC}_{-1} \left(\pi\left(\mc{L}(-1, 1)\otimes \mr{MC}_{-1}(\Sym^2 \mc{G}) \right)\right)\right)
\]
has the desired local monodromies, where $\mc{G}$ is the $\mr{GL}_2$-rigid sheaf (a well-known classical example) having local monodromies
\[
U(2) \qquad -U(2) \qquad U(2).
\]
Here too the $\mr{Spin}_7$-lifts are $\mr{Spin}_7$-rigid but not $\mr{SO}_8$-rigid, reflecting the increased difficulty of showing they are motivic. Note too that even were we to find a $\mr{Spin}_7$-local system $\tilde{\rho}$ whose local monodromies in the 8-dimensional spin representation were the desired ones, the argument as in Theorem \ref{main} would not quite work: the original $\mr{SO}_7$-local system is not \textit{a priori} physically rigid, which means that the comparison of local monodromies does not suffice to test global isomorphism of this original local system and the $\mr{SO}_7$-reduction of $\tilde{\rho}$. I expect that this difficulty too can be handled.
\end{eg}
\bibliographystyle{amsalpha}
\bibliography{biblio.bib}

\end{document}